\documentclass[11pt]{article}
\usepackage{amssymb}
\newtheorem{lemma}{Lemma}[section]
\newtheorem{definition}[lemma]{Definition}

\newtheorem{theorem}[lemma]{Theorem}
\newtheorem{proposition}[lemma]{Proposition}

\newtheorem{example}[lemma]{Example}

\newtheorem{note}[lemma]{Remark}

\def\endproof{\hfill$\Box$}

\def\endproof{\hfill$\Box$}

\title{Distributions of countable models  \\ of disjoint unions of Ehrenfeucht theories\footnote{This research was partially
supported by Russian Foundation for Basic Researches (Project No.
17-01-00531-a).}}
\author{S.V. Sudoplatov}
\date{}

\begin{document}

\maketitle

\begin{abstract}
We describe Rudin--Keisler preorders and distribution functions of
numbers of limit models for disjoint unions of Ehrenfeucht
theories. Decomposition formulas for these distributions are
found.

\end{abstract}

In \cite{KulSudRK}, a description is obtained for distributions of
countable models of quite $o$-minimal Ehrenfeucht theories in
terms of Rudin--Keisler preorders and distribution functions of
numbers of limit models. In the present paper, using a general
theory of classification of countable models of complete theories
\cite{SuCCMCT1, SuCCMCT2} as well as the description
\cite{{KulSudRK}} of specificity for quite $o$-minimal theories,
we describe distributions of countable models of disjoint unions
of Ehrenfeucht theories in terms of Rudin--Keisler preorders and
distribution functions of numbers of limit models. Besides, we
derive decomposition formulas for these distributions.

\section{Preliminaries}

In this section we give the necessary information from
\cite{SuCCMCT1, SuCCMCT2}.

Recall that the number of pairwise non-isomorphic models of theory
$T$ and of cardinality $\lambda$ is denoted by
$I(T,\lambda)$\index{$I(T,\lambda)$}.

\begin{definition} {\rm \cite{Mi2} A theory $T$ is called
\emph{Ehrenfeucht} if $1<I(T,\omega)<\omega$.}
\end{definition}

\begin{definition} {\rm \cite{Be}  type $p(\bar{x})\in S(T)$ is said to be
\emph{powerful}\index{Type!powerful} in a theory $T$ if every
model ${\cal M}$ of $T$ realizing $p$ also realizes every type
$q\in S(T)$, that is, ${\cal M}\models S(T)$.}
\end{definition}

Since for any type $p\in S(T)$ there exists a countable model
${\cal M}$ of $T$, realizing $p$, and the model ${\cal M}$
realizes exactly countably many types, the availability of a
powerful type implies that $T$ is
\emph{small}\index{Theory!small}, that is, the set $S(T)$ is
countable. Hence for any type $q\in S(T)$ and its realization
$\bar{a}$, there exists a {\em prime model} ${\cal
M}({\bar{a}})$\index{${\cal M}({\bar{a}})$} {\em over}
$\bar{a}$,\index{Model!prime over a tuple} i.~e., a model of $T$
containing $\bar{a}$ with ${\cal M}(\bar{a})\models q(\bar{a})$
and such that ${\cal M}({\bar{a}})$ is elementarily embeddable to
any model realizing the type $q$. Since all prime models over
realizations of $q$ are isomorphic, we denote these models by
${\cal M}_q$\index{${\cal M}_q$}. Models ${\cal M}_q$ are called
{\em almost prime}\index{Model!almost prime} or {\em
$q$-prime}.\index{Model!$q$-prime}

\begin{definition} {\rm \cite{SuCCMCT1, Lr3, Su041} Let $p$ and $q$ be types in $S(T)$. We say that the type $p$
\emph{is dominated by a type}\index{Type!dominated} $q$, or $p$
\emph{does not exceed $q$ under the Rudin--Keisler
preorder}\index{Type!not exceed $q$ under the Rudin--Keisler
preorder}\index{Preorder!Rudin--Keisler} (written $p\leq_{\rm RK}
q$\index{$p\leq_{\rm RK}q$}\index{$\leq_{\rm RK}$}), if ${\cal
M}_q\models p$, that is, ${\cal M}_p$ is an elementary submodel of
${\cal M}_q$ (written ${\cal M}_p\preceq{\cal M}_q$). Besides, we
say that a~model ${\cal M}_p$ \emph{is dominated by a
model}\index{Model!dominated} ${\cal M}_q$, or ${\cal M}_p$
\emph{does not exceed ${\cal M}_q$ under the Rudin--Keisler
preorder}\index{Model!not exceed $q$ under the Rudin--Keisler
preorder}, and write ${\cal M}_p\leq_{\rm RK}{\cal
M}_q$\index{${\cal M}_p\leq_{\rm RK}{\cal M}_q$}.}
\end{definition}

Syntactically,  the  condition \ $p\leq_{\rm RK}q$  (and  hence
also ${\cal M}_p\leq_{\rm RK}{\cal M}_q$)  is expressed thus:
there exists a formula $\varphi(\bar{x},\bar{y})$ such that the
set $q(\bar{y})\cup\{\varphi(\bar{x},\bar{y})\}$ is consistent and
$q(\bar{y})\cup\{\varphi(\bar{x},\bar{y})\}\vdash p(\bar{x})$.
Since we deal with a small theory (there are only countably many
types over any tuple $\bar{a}$ and so any consistent formula with
parameters in $\bar{a}$ is deducible from a principal formula with
parameters in $\bar{a}$), $\varphi(\bar{x},\bar{y})$ can be chosen
so that for any formula $\psi(\bar{x},\bar{y})$, the set
$q(\bar{y})\cup\{\varphi(\bar{x},\bar{y}),\psi(\bar{x},\bar{y})\}$
being consistent implies that
$q(\bar{y})\cup\{\varphi(\bar{x},\bar{y})\}\vdash\psi(\bar{x},\bar{y})$.
In this event the formula $\varphi(\bar{x},\bar{y})$ is said to be
\emph{$(q,p)$-principal}\index{Formula!$(q,p)$-principal}.

\medskip
\begin{definition} {\rm \cite{SuCCMCT1, Lr3, Su041} Types $p$ and $q$ are said to be
\emph{domination-equivalent}\index{Types!domination-equivalent},
\emph{realization-equivalent}\index{Types!realization-equivalent},
\emph{Rudin--Keisler equivalent}\index{Types!Rudin--Keisler
equivalent}, or \ \emph{${\rm RK}$-equivalent}\index{Types!${\rm
RK}$-equivalent} \ (written \ $p\sim_{\rm RK} q$\index{$p\sim_{\rm
RK} q$}) if~$p\leq_{\rm RK} q$ and $q\leq_{\rm RK} p$. Models
${\cal M}_p$ and ${\cal M}_q$ are said to be
\emph{domination-equivalent}\index{Models!domination-equivalent},
\ \emph{Rudin--Keisler \ equivalent}\index{Models!Rudin--Keisler
equivalent}, \ or \ \emph{${\rm
RK}$-equivalent}\index{Models!${\rm RK}$-equivalent} \ (written \
\ ${\cal M}_p\sim_{\rm RK}{\cal M}_q$\index{${\cal M}_p\sim_{\rm
RK}{\cal M}_q$}).

As in \cite{Ta4}, types $p$ and $q$ are said to be \emph{strongly
domination-equivalent}\index{Types!domination-equivalent!strongly},
\emph{strongly
realization-equivalent}\index{Types!realization-equivalent!strongly},
\emph{strongly Rudin--Keisler
equivalent}\index{Types!Rudin--Keisler equivalent!strongly}, or
\emph{strongly ${\rm RK}$-equivalent}\index{Types!${\rm
RK}$-equivalent!strongly} (written \ $p\equiv_{\rm
RK}q$\index{$p\equiv_{\rm RK} q$}) if for some realizations
$\bar{a}$ and $\bar{b}$ of $p$ and $q$ respectively, both ${\rm
tp}(\bar{b}/\bar{a})$ and ${\rm tp}(\bar{a}/\bar{b})$ are
principal. Models ${\cal M}_p$ and ${\cal M}_q$ are said to be
\emph{strongly
domination-equivalent}\index{Models!domination-equivalent!strongly},
\emph{strongly Rudin--Keisler
equivalent}\index{Models!Rudin--Keisler equivalent}, or
\emph{strongly ${\rm RK}$-equivalent}\index{Models!${\rm
RK}$-equivalent!strongly} (written ${\cal M}_p\equiv_{\rm RK}{\cal
M}_q$\index{${\cal M}_p\equiv_{\rm RK}{\cal M}_q$}).}
\end{definition}

\medskip
Clearly, domination relations form preorders, and (strong)
do\-mi\-na\-tion-equi\-va\-lence \ relations \ are \ equivalence \
relations. \ Here, \ ${\cal M}_p\equiv_{\rm RK}{\cal M}_q$ implies
${\cal M}_p\sim_{\rm RK}{\cal M}_q$.

If ${\cal M}_p$ and ${\cal M}_q$ are not domina\-tion-equivalent
then they are non-isomorphic. Moreover, non-isomorphic models may
be found among domination-equivalent ones.

In \ Ehrenfeucht \ examples, \ models \ ${\cal
M}^n_{p_0},\ldots,{\cal M}^n_{p_{n-3}}$ \ are
domina\-tion-equivalent but pairwise non-isomorphic.

A syntactic characterization for the model isomorphism between
${\cal M}_p$ and ${\cal M}_q$ is given by the following
proposition. It asserts that the existence of an isomorphism
between ${\cal M}_p$ and ${\cal M}_q$ is equivalent to the strong
domination-equivalence of these models.

\medskip
%\begin{statement}\label{st117_1117}
\begin{proposition} {\rm \cite{SuCCMCT1, Su041}} For any types
$p(\bar{x})$ and $q(\bar{y})$ of a small theory $T$, the following
conditions are equivalent:

{\rm (1)} the models ${\cal M}_p$ and ${\cal M}_q$ are isomorphic;

{\rm (2)} the models ${\cal M}_p$ and ${\cal M}_q$ are strongly
domination-equivalent;

{\rm (3)} there exist $(p,q)$- and $(q,p)$-principal formulas
$\varphi_{p,q}(\bar{y},\bar{x})$ and
$\varphi_{q,p}(\bar{x},\bar{y})$ respectively, such that the set
$$
p(\bar{x})\cup q(\bar{y})\cup\{\varphi_{p,q}(\bar{y},\bar{x}),
\varphi_{q,p}(\bar{x},\bar{y})\}$$ is consistent;

{\rm (4)} there exists a $(p,q)$- and $(q,p)$-principal formula
$\varphi(\bar{x},\bar{y})$, such that the set
$$
p(\bar{x})\cup q(\bar{y})\cup\{\varphi(\bar{x},\bar{y})\}$$ is
consistent.
\end{proposition}

\begin{definition} {\rm \cite{SuCCMCT1, Su041}  Denote  by \ ${\rm RK}(T)$\index{${\rm RK}(T)$}  the  set
${\bf PM}$\index{${\bf PM}$}  of  isomorphism  types  of models
${\cal M}_p$, $p\in S(T)$, on which the relation of domination is
induced by $\leq_{\rm RK}$, a relation deciding domination among
${\cal M}_p$, that is, ${\rm RK}(T)=\langle{\bf PM};\leq_{\rm
RK}\rangle$. We say that isomorphism types ${\bf M}_1,{\bf
M}_2\in{\bf PM}$ are \emph{domination-equivalent} (written ${\bf
M}_1\sim_{\rm RK}{\bf M}_2$\index{${\bf M}_1\sim_{\rm RK}{\bf
M}_2$}) if so are their representatives.}
\end{definition}

Clearly, the preordered set ${\rm RK}(T)$ has a least element,
which is an isomorphism type of a prime model.

\medskip
%\begin{statement}\label{st118_1119}
\begin{proposition} {\rm \cite{SuCCMCT1, Su041}} If $I(T,\omega)<\omega$ then ${\rm
RK}(T)$ is a finite preordered set whose factor set ${\rm
RK}(T)/\!\!\!\sim_{\rm RK}$, with respect to
do\-mi\-na\-tion-equiva\-lence $\sim_{\rm RK}$, forms a partially
ordered set with a greatest element.
\end{proposition}

\begin{definition} {\rm \cite{SuCCMCT1, SuCCMCT2, Su041, Su08}
A model ${\cal M}$ of a theory $T$ is called {\em
limit}\index{Model!limit} if ${\cal M}$ is not prime over tuples
and ${\cal M}=\bigcup\limits_{n\in\omega}{\cal M}_n$ for some
elementary chain $({\cal M}_n)_{n\in\omega}$ of prime models of
$T$ over tuples. In this case the model ${\cal M}$ is said to be
{\em limit over a sequence ${\bf q}$ of
types}\index{Model!limit!over a sequence of types} or  {\em ${\bf
q}$-limit}\index{Model!${\bf q}$-limit}, where ${\bf
q}=(q_n)_{n\in\omega}$, ${\cal M}_n={\cal M}_{q_n}$, $n\in\omega$.
If the sequence ${\bf q}$ contains unique type $q$ then the ${\bf
q}$-limit model is called {\em limit over the type $q$}.}
\end{definition}

Denote \ by \ $I_p(T,\omega)$\index{$I_p(T,\omega)$} \ the \
number \ of \ pairwise \ non-isomorphic \ countable models of the
theory $T$, each of which is prime over a tuple, by $I_l(T)$ the
number of limit models of $T$, and by $I_l(T,q)$ the number of
limit models over a type $q\in S(T)$.

\begin{definition}
{\rm \cite{SuCCMCT2, Su08} A  theory  $T$  is  called  {\em
$p$-categorical}\index{Theory!$p$-categorical}  (respectively,
{\em $l$-categorical}\index{Theory!$l$-categorical},  {\em
$p$-Ehren\-feucht}\index{Theory!$p$-Ehrenfeucht}, and {\em
$l$-Ehrenfeucht})\index{Theory!$l$-Ehrenfeucht}  if
$I_p(T,\omega)=1$ \ (respectively, \ $I_l(T)=1$,
$1<I_p(T,\omega)<\omega$, $1<I_l(T)<\omega$).}
\end{definition}

Clearly, a small theory $T$ is $p$-categorical if and only if $T$
countably categorical, and if and only if $I_l(T)=0$; $T$ is
$p$-Ehrenfeucht if and only if the structure ${\rm RK}(T)$ finite
and has at least two elements; and $T$ is $p$-Ehrenfeucht with
$I_l(T)<\omega$ if and only if $T$ is Ehrenfeucht.

\medskip
Let $\widetilde{\bf M}\in{\rm RK}(T)/\!\!\sim_{\rm RK}$ be the
class consisting of isomorphism types of domination-equivalent
models ${\cal M}_{p_1},\ldots,{\cal M}_{p_n}$. Denote by ${\rm
IL}(\widetilde{\bf M})$ \index{${\rm IL}(\widetilde{\bf M})$} the
number of equivalence classes of models each of which is limit
over some type $p_i$.

\begin{theorem}\label{th1121_1136} {\rm \cite{SuCCMCT1, Su041}}
For any countable complete theory $T$, the following conditions
are equivalent:

{\rm (1)} $I(T,\omega)<\omega$;

{\rm (2)} $T$  is  small,  $|{\rm RK}(T)|<\omega$  and  ${\rm
IL}(\widetilde{\bf M})<\omega$  for  any  $\widetilde{\bf
M}\in{\rm RK}(T)/\!\!\sim_{\rm RK}$.

If $(1)$ or $(2)$ holds then $T$ possesses the following
properties:

{\rm (a)} ${\rm RK}(T)$ has a least element ${\bf M}_0$ {\rm (}an
isomorphism type of a prime model{\rm )} and ${\rm
IL}(\widetilde{{\bf M}_0})=0$;

{\rm (b)} ${\rm RK}(T)$ has a greatest $\sim_{\rm RK}$-class
$\widetilde{{\bf M}_1}$ {\rm (}a class of isomorphism types of all
prime models over realizations of powerful types{\rm )} and $|{\rm
RK}(T)|>1$ implies ${\rm IL}(\widetilde{{\bf M}_1})\geq 1$;

{\rm (c)} if $|\widetilde{\bf M}|>1$ then ${\rm IL}(\widetilde{\bf
M})\geq 1$.

Moreover, the following {\sl decomposition
formula}\index{Formula!decomposition} holds:
\begin{equation}\label{eqmain}
I(T,\omega)=|{\rm RK}(T)|+\sum_{i=0}^{|{\rm RK}(T)/\sim_{\rm
RK}|-1} {\rm IL}(\widetilde{{\bf M}_i}),
\end{equation}
where $\widetilde{{\bf M}_0},\ldots, \widetilde{{\bf M}_{|{\rm
RK}(T)/\sim_{\rm RK}|-1}}$ are all elements of the partially
ordered set ${\rm RK}(T)/\!\!\sim_{\rm RK}$.
\end{theorem}

In Figure \ref{fig1}, {\it a} and {\it b}, possible variants for
Hasse diagrams of~Rudin--Keisler preorders $\leq_{\rm RK}$ and
values of distribution functions ${\rm IL}$ of numbers of limit
models on $\sim_{\rm RK}$-equivalence classes are represented for
the cases $I(T,\omega)=3$ and $I(T,\omega)=4$. In
Figure~\ref{fig2}, corresponding configurations for
$I(T,\omega)=5$ are shown.

\begin{figure}
\vspace{5mm}
\begin{center}
\unitlength 4mm
\begin{picture}(22,6.5)(0,1.0)
{\footnotesize\put(2,2.5){\line(0,5){5}}
\put(2,2.5){\makebox(0,0)[cc]{$\bullet$}}
\put(2,7.5){\makebox(0,0)[cc]{$\bullet$}}
\put(2,2.5){\circle{0.6}} \put(2,7.5){\circle{0.6}}
\put(3,2.5){\makebox(0,0)[cr]{$0$}}
\put(3,7.5){\makebox(0,0)[cr]{$1$}}

\put(9,2.5){\line(0,5){5}}
\put(9,2.5){\makebox(0,0)[cc]{$\bullet$}}
\put(9,7.5){\makebox(0,0)[cc]{$\bullet$}}
\put(9,2.5){\circle{0.6}} \put(9,7.5){\circle{0.6}}
\put(10,2.5){\makebox(0,0)[cr]{$0$}}
\put(10,7.5){\makebox(0,0)[cr]{$2$}}

\put(14.5,2.5){\line(-2,5){2}} \put(14.5,2.5){\line(2,5){2}}
\put(14.5,2.5){\makebox(0,0)[cc]{$\bullet$}}
\put(12.5,7.5){\makebox(0,0)[cc]{$\bullet$}}
\put(16.5,7.5){\makebox(0,0)[cc]{$\bullet$}}
\put(14.5,2.5){\circle{0.6}} \put(14.5,7.5){\oval(5,1)}
\put(15.5,2.5){\makebox(0,0)[cr]{$0$}}
\put(17.7,7.5){\makebox(0,0)[cr]{$1$}}

\put(20,2.5){\line(0,5){5}}
\put(20,2.5){\makebox(0,0)[cc]{$\bullet$}}
\put(20,5){\makebox(0,0)[cc]{$\bullet$}}
\put(20,7.5){\makebox(0,0)[cc]{$\bullet$}}
\put(20,2.5){\circle{0.6}}
\put(20,5){\circle{0.6}}\put(20,7.5){\circle{0.6}}
\put(21,2.5){\makebox(0,0)[cr]{$0$}}
\put(21,5){\makebox(0,0)[cr]{$0$}}
\put(21,7.5){\makebox(0,0)[cr]{$1$}}

\put(2,0){\makebox(0,0)[bc]{\textit{a}}}
\put(14.5,0){\makebox(0,0)[bc]{\textit{b}}}}
\end{picture}
\end{center}
\caption{} \label{fig1}
\end{figure}

\begin{figure}[t]
\begin{center}
\unitlength 5mm
\begin{picture}(22,13)(0,2.5)
{\footnotesize\put(2,10){\line(0,5){5}}
\put(2,10){\makebox(0,0)[cc]{$\bullet$}}
\put(2,15){\makebox(0,0)[cc]{$\bullet$}} \put(2,10){\circle{0.6}}
\put(2,15){\circle{0.6}} \put(3,10){\makebox(0,0)[cr]{$0$}}
\put(3,15){\makebox(0,0)[cr]{$3$}}

\put(8,10){\line(0,5){5}} \put(8,10){\makebox(0,0)[cc]{$\bullet$}}
\put(8,12.5){\makebox(0,0)[cc]{$\bullet$}}
\put(8,15){\makebox(0,0)[cc]{$\bullet$}} \put(8,10){\circle{0.6}}
\put(8,12.5){\circle{0.6}} \put(8,15){\circle{0.6}}
\put(9,10){\makebox(0,0)[cr]{$0$}}
\put(9,12.5){\makebox(0,0)[cr]{$1$}}
\put(9,15){\makebox(0,0)[cr]{$1$}}

\put(14,10){\line(0,5){5}}
\put(14,10){\makebox(0,0)[cc]{$\bullet$}}
\put(14,12.5){\makebox(0,0)[cc]{$\bullet$}}
\put(14,15){\makebox(0,0)[cc]{$\bullet$}}
\put(14,10){\circle{0.6}} \put(14,12.5){\circle{0.6}}
\put(14,15){\circle{0.6}} \put(15,10){\makebox(0,0)[cr]{$0$}}
\put(15,12.5){\makebox(0,0)[cr]{$0$}}
\put(15,15){\makebox(0,0)[cr]{$2$}}

\put(20,10){\line(0,5){5}}
\put(20,10){\makebox(0,0)[cc]{$\bullet$}}
\put(20,11.6666){\makebox(0,0)[cc]{$\bullet$}}
\put(20,13.3333){\makebox(0,0)[cc]{$\bullet$}}
\put(20,15){\makebox(0,0)[cc]{$\bullet$}}
\put(20,10){\circle{0.6}}
\put(20,11.6666){\circle{0.6}}\put(20,13.3333){\circle{0.6}}
\put(20,15){\circle{0.6}} \put(21,10){\makebox(0,0)[cr]{$0$}}
\put(21,11.6666){\makebox(0,0)[cr]{$0$}}
\put(21,13.3333){\makebox(0,0)[cr]{$0$}}
\put(21,15){\makebox(0,0)[cr]{$1$}}

\put(2,2.5){\line(-1,3){1.5}} \put(2,2.5){\line(1,3){1.5}}
\put(2,2.5){\makebox(0,0)[cc]{$\bullet$}}
\put(0.5,7){\makebox(0,0)[cc]{$\bullet$}}
\put(3.5,7){\makebox(0,0)[cc]{$\bullet$}}
\put(2,2.5){\circle{0.6}} \put(2,7){\oval(4,1)}
\put(3,2.5){\makebox(0,0)[cr]{$0$}}
\put(4.5,7){\makebox(0,0)[cr]{$2$}}

\put(8,2.5){\line(-1,3){1.5}} \put(8,2.5){\line(1,3){1.5}}
\put(8,2.5){\line(0,1){4.5}}
\put(8,2.5){\makebox(0,0)[cc]{$\bullet$}}
\put(6.5,7){\makebox(0,0)[cc]{$\bullet$}}
\put(9.5,7){\makebox(0,0)[cc]{$\bullet$}}
\put(8,7){\makebox(0,0)[cc]{$\bullet$}} \put(8,2.5){\circle{0.6}}
\put(8,7){\oval(4,1)} \put(9,2.5){\makebox(0,0)[cr]{$0$}}
\put(10.5,7){\makebox(0,0)[cr]{$1$}}

\put(14,2.5){\line(0,1){3}} \put(14,5.5){\line(1,1){1.5}}
\put(14,5.5){\line(-1,1){1.5}}
\put(14,2.5){\makebox(0,0)[cc]{$\bullet$}}
\put(12.5,7){\makebox(0,0)[cc]{$\bullet$}}
\put(15.5,7){\makebox(0,0)[cc]{$\bullet$}}
\put(14,5.5){\makebox(0,0)[cc]{$\bullet$}}
\put(14,2.5){\circle{0.6}} \put(14,5.5){\circle{0.6}}
\put(14,7){\oval(4,1)} \put(15,2.5){\makebox(0,0)[cr]{$0$}}
\put(15,5.5){\makebox(0,0)[cr]{$0$}}
\put(16.5,7){\makebox(0,0)[cr]{$1$}}

\put(20,2.5){\line(-3,5){1.5}} \put(20,2.5){\line(3,5){1.5}}
\put(18.5,5){\line(3,5){1.5}} \put(21.5,5){\line(-3,5){1.5}}
\put(20,2.5){\makebox(0,0)[cc]{$\bullet$}}
\put(18.5,5){\makebox(0,0)[cc]{$\bullet$}}
\put(21.5,5){\makebox(0,0)[cc]{$\bullet$}}
\put(20,7.5){\makebox(0,0)[cc]{$\bullet$}}
\put(20,2.5){\circle{0.6}} \put(20,7.5){\circle{0.6}}
\put(18.5,5){\circle{0.6}} \put(21.5,5){\circle{0.6}}
\put(21,2.5){\makebox(0,0)[cr]{$0$}}
\put(21,7.5){\makebox(0,0)[cr]{$1$}}
\put(19.3,5){\makebox(0,0)[cr]{$0$}}
\put(21,5){\makebox(0,0)[cr]{$0$}}

}\end{picture}
\end{center}
\caption{} \label{fig2}
\end{figure}

\begin{definition} {\rm \cite{Wo} The {\em disjoint union}\index{Disjoint union!of
structures} $\bigsqcup\limits_{n\in\omega}{\cal
M}_n$\index{$\bigsqcup\limits_{n\in\omega}{\cal M}_n$} of pairwise
disjoint structures ${\cal M}_n$ for pairwise disjoint predicate
languages $\Sigma_n$, $n\in\omega$, is the structure of language
$\bigcup\limits_{n\in\omega}\Sigma_n\cup\{P^{(1)}_n\mid
n\in\omega\}$ with the universe $\bigsqcup\limits_{n\in\omega}
M_n$, $P_n=M_n$, and interpretations of predicate symbols in
$\Sigma_n$ coinciding with their interpretations in ${\cal M}_n$,
$n\in\omega$. The {\em disjoint union of theories}\index{Disjoint
union!of theories} $T_n$ for pairwise disjoint languages
$\Sigma_n$ accordingly, $n\in\omega$, is the theory
$$\bigsqcup\limits_{n\in\omega}T_n\rightleftharpoons{\rm Th}\left(\bigsqcup\limits_{n\in\omega}{\cal M}_n\right),$$
where\index{$\bigsqcup\limits_{n\in\omega}T_n$} ${\cal M}_n\models
T_n$, $n\in\omega$.}
\end{definition}

Clearly, the theory $T_1\sqcup T_2$ does not depend on choice of
disjoint union ${\cal M}_1\sqcup {\cal M}_2$ of models ${\cal
M}_1\models T_1$ and ${\cal M}_2\models T_2$. Besides, the
cardinality of ${\rm RK}(T_1\sqcup T_2)$ is equal to the product
of cardinalities for ${\rm RK}(T_1)$ and ${\rm RK}(T_2)$, and the
relation $\leq_{\rm RK}$ on ${\rm RK}(T_1\sqcup T_2)$ equals the
Pareto relation \cite{SO1} defined by preorders in  ${\rm
RK}(T_1)$ and ${\rm RK}(T_2)$. Indeed, each type $p(\bar{x})$ of
$T_1\sqcup T_2$ is isolated by set consisting of some types
$p_1(\bar{x}^1)$ and $p_2(\bar{x}^2)$ of theories $T_1$ and $T_2$
respectively, as well as of formulas $P^1(x^1_i)$ and $P^2(x^2_j)$
for all coordinates in tuples $\bar{x}^1$ and $\bar{x}^2$. For
types $p(\bar{x})$ and $p'(\bar{y})$ of $T_1\sqcup T_2$, we have
$p(\bar{x})\leq_{\rm RK} p'(\bar{y})$ if and only if
$p_1(\bar{x}^1)\leq_{\rm RK} p'_1(\bar{y}^1)$ (in $T_1$) and
$p_2(\bar{x}^2)\leq_{\rm RK} p'_2(\bar{y}^2)$ (in $T_2$).

Thus, the following proposition holds.

\medskip
%\begin{statement}
\begin{proposition}\label{st676} {\rm \cite{SuCCMCT2, Su09Irk}}
For any small theories $T_1$ and $T_2$ of disjoint predicate
languages $\Sigma_1$ and $\Sigma_2$ respectively, the theory
$T_1\sqcup T_2$ is mutually ${\rm RK}$-coordinated with respect to
its restrictions to $\Sigma_1$ and $\Sigma_2$. The cardinality of
${\rm RK}(T_1\sqcup T_2)$ is equal to the product of cardinalities
for ${\rm RK}(T_1)$ and ${\rm RK}(T_2)$, i.~e.,
\begin{equation}\label{s64} I_p(T_1\sqcup
T_2,\omega)=I_p(T_1,\omega)\cdot I_p(T_2,\omega),
\end{equation} and the relation
$\leq_{\rm RK}$ on ${\rm RK}(T_1\sqcup T_2)$ equals the Pareto
relation defined by preorders in ${\rm RK}(T_1)$ and ${\rm
RK}(T_2)$.
\end{proposition}

\begin{note}\label{no6742}
{\rm \cite{SuCCMCT2, Su09Irk} An isomorphism \ of limit \ models \
of theory \ $T_1\sqcup T_2$ is defined by isomorphisms of
restrictions of these models to the sets $P_1$ and $P_2$. In this
case, a countable model is limit if and only if some its
restriction (to $P_1$ or to $P_2$) is limit and the following
equality holds:
\begin{equation}\label{s65}
I(T_1\sqcup T_2,\omega)=I(T_1,\omega)\cdot I(T_2,\omega).
\end{equation}
Thus, the operation $\sqcup$ preserves both $p$-Ehrenfeuchtness
and $l$-Ehrenfeuchtness (if components are $p$-Ehrenfeucht), and,
by $(\ref{s65})$, we obtain the equality
\begin{equation}\label{s66}
I_l(T_1\sqcup T_2)=I_l(T_1)\cdot I_p(T_2,\omega)+
I_p(T_1,\omega)\cdot I_l(T_2)+I_l(T_1)\cdot I_l(T_2).
\end{equation}}
\end{note}

\section{Distributions of countable models}

In this section, using Theorem \ref{th1121_1136} and Proposition
\ref{st676}, we give a description of Rudin--Keisler preorders and
distribution functions of numbers of limit models for disjoint
unions  $T_1\sqcup T_2$ of Ehrenfeucht theories $T_1$ and $T_2$,
as well as propose representations of this distributions, based on
the decomposition formula (\ref{eqmain}).

Using the formulas (\ref{eqmain})--(\ref{s66}) we obtain the
following equalities:
$$
I(T_1,\omega)\cdot I(T_2,\omega)=I(T_1\sqcup
T_2,\omega)=I_p(T_1\sqcup T_2,\omega)+I_l(T_1\sqcup T_2)=$$
\begin{equation}\label{s67}
=I_p(T_1,\omega)\cdot I_p(T_2,\omega)+I_l(T_1\sqcup
T_2)=I_l(T_1)\cdot I_p(T_2,\omega)+ I_p(T_1,\omega)\cdot
I_l(T_2)+I_l(T_1)\cdot I_l(T_2)
\end{equation}
implying
\begin{equation}\label{s68}
I(T_1,\omega)\cdot I(T_2,\omega)= I_p(T_1,\omega)\cdot
I_p(T_2,\omega)+I_l(T_1)\cdot I_p(T_2,\omega)+
I_p(T_1,\omega)\cdot I_l(T_2)+I_l(T_1)\cdot I_l(T_2).
\end{equation}

In view of Proposition \ref{st676} the Hasse diagrams for
distributions of countable models for disjoint unions $T_1\sqcup
T_2$ of Ehrenfeucht theories $T_1$ and $T_2$ are constructed as
images of Pareto relations for these theories $T_1$ and $T_2$.
Here, $\sim_{\rm RK}$-equivalent vertices for $T_1$ and $T_2$ are
transformed to $\sim_{\rm RK}$-equivalent pairs for $T_1\sqcup
T_2$. Hence, each $\sim_{\rm RK}$-class for $T_1$, consisting of
$k$ vertices, united with a $\sim_{\rm RK}$-class for $T_2$,
consisting of $m$ vertices, produces a $\sim_{\rm RK}$-class
$\tilde{z}$ for $T_1\sqcup T_2$, consisting of $km$ vertices.
Thus, in the formula (\ref{s68}), the value $I_p(T_1,\omega)\cdot
I_p(T_2,\omega)$ is represented as a sum of products
$|\tilde{x}|\cdot|\tilde{y}|$ for each equivalence class
$\tilde{x}$ in ${\rm RK}(T_1)$ and each equivalence class
$\tilde{y}$ in ${\rm RK}(T_2)$:
\begin{equation}\label{s69}
I_p(T_1\sqcup T_2,\omega)= I_p(T_1,\omega)\cdot
I_p(T_2,\omega)=\sum\limits_{\tilde{x}\in {\rm RK}(T_1),
\tilde{y}\in {\rm RK}(T_2)}|\tilde{x}|\cdot|\tilde{y}|.
\end{equation}

Following the formula (\ref{s66}), each $\sim_{\rm RK}$-class
$\tilde{z}$ has some number $I_l(\tilde{z})$ of limit models over
types defining that class. This number is expressed by the numbers
$I_l(\tilde{x})$ and $I_l(\tilde{y})$ of limit models for the
$\sim_{\rm RK}$-class $\tilde{x}$ in ${\rm RK}(T_1)$ and the
$\sim_{\rm RK}$-class $\tilde{y}$ in ${\rm RK}(T_2)$, generating
$\tilde{z}$, by the following formula:
\begin{equation}\label{s70}
I_l(\tilde{z})=I_l(\tilde{x})\cdot |\tilde{y}|+ |\tilde{x}|\cdot
I_l(\tilde{y})+I_l(\tilde{x})\cdot I_l(\tilde{y}).
\end{equation}

By (\ref{s69}) and (\ref{s70}), for the theory $T_1\sqcup T_2$,
the decomposition formula (\ref{eqmain}) has the following form:
\begin{equation}\label{s71}
I(T_1\sqcup T_2,\omega)= \sum\limits_{\tilde{x},
\tilde{y}}|\tilde{x}|\cdot|\tilde{y}|+\sum\limits_{\tilde{x},
\tilde{y}}\left(I_l(\tilde{x})\cdot |\tilde{y}|+ |\tilde{x}|\cdot
I_l(\tilde{y})+I_l(\tilde{x})\cdot I_l(\tilde{y})\right).
\end{equation}

Notice that the graph $\Gamma$ for the Pareto relation
correspondent to the Rudin--Keisler preorder of the theory
$T_1\sqcup T_2$ is represented as the product of the graphs
$\Gamma_1$ and $\Gamma_2$ for the Rudin--Keisler preorders of the
theories $T_1$ and $T_2$, and $\Gamma_1\times\Gamma_2$ is a
(Boolean) lattice if and only if $\Gamma_1$ and $\Gamma_2$ are
(Boolean) lattices.

Thus, the following theorem holds, generalizing Theorem 24 in
\cite{KulSudRK}.

\begin{theorem}\label{thchar2}
For any Ehrenfeucht theories $T_1$ and $T_2$ with graphs
$\Gamma_1$ and $\Gamma_2$ of Rudin--Keisler preorders the theory
$T_1\sqcup T_2$ has the Rudin--Keisler preorder, represented by
the product $\Gamma_1\times\Gamma_2$, and the decomposition
formula of the form {\rm (\ref{s71})}. The structure
$\Gamma_1\times\Gamma_2$ for $T_1\sqcup T_2$ forms a {\rm
(}Boolean{\rm )} lattice if and only if $\Gamma_1$ and $\Gamma_2$
form {\rm (}Boolean{\rm )} lattices.
\end{theorem}

The following proposition shows that the function $I_l(\cdot)$ is
monotone with respect to disjoint unions of Ehrenfeucht theories.

\begin{proposition}\label{prchar2}
The functions $I_l(\tilde{x})$ and $I_l(\tilde{y})$ of numbers of
limit models for $\sim_{\rm RK}$-classes $\tilde{x}$ in ${\rm
RK}(T_1)$ and $\tilde{y}$ in ${\rm RK}(T_2)$ monotonically
increase {\rm (}do not decrease{\rm )} with respect to
Rudin--Keisler preorders, having monotonically increasing {\rm
(}non-decreasing{\rm )} cardinalities $|\tilde{x}|$ and
$|\tilde{y}|$, if and only if the function $I_l(\tilde{z})$ of
numbers of limit models for $\sim_{\rm RK}$-classes $\tilde{z}$ in
${\rm RK}(T_1\sqcup T_2)$ monotonically increase {\rm (}do not
decrease{\rm )} with respect to Rudin--Keisler preorder, having
monotonically increasing {\rm (}non-decreasing{\rm )}
cardinalities $|\tilde{z}|$.
\end{proposition}

Proof. Assume that the cardinalities $|\tilde{x}|$ and
$|\tilde{y}|$ monotonically increase {\rm (}do not decrease{\rm )}
with respect to Rudin--Keisler preorders. If the functions
$I_l(\tilde{x})$ and $I_l(\tilde{y})$ monotonically increase {\rm
(}do not decrease{\rm )} with respect to Rudin--Keisler preorders
and $\tilde{z}_1<_{\rm RK}\tilde{z}_2$ ($\tilde{z}_1\leq_{\rm
RK}\tilde{z}_2$), then $I_l(\tilde{z}_1)<I_l(\tilde{z}_2)$
($I_l(\tilde{z}_1)\leq I_l(\tilde{z}_2)$) in view of the formula
(\ref{s70}).

The reverse implication takes place, since ${\rm RK}(T_1)$ and
${\rm RK}(T_2)$ are isomorphic to substructures of the structure
${\rm RK}(T_1\sqcup T_2)$.
\endproof

\medskip
In addition to the examples of Hasse diagrams given in
\cite{KulSudRK}, we have a series of new examples. Below we
describe some of them.

\begin{example}\label{exx1}\rm
Consider the disjoint union of theory $T_1$ with the Hasse diagram
shown in Fig.~\ref{fig1}, {\it a} and of theory $T_2$ with the
first Hasse diagram shown in Fig.~\ref{fig1}, {\it b}. By Theorem
\ref{thchar2} we have the theory $T_1\sqcup T_2$ with $3\cdot
4=12$ countable models, having the Boolean lattice with $2\cdot
2=4$ prime models over finite sets and with 8 limit models. The
decomposition formula (\ref{s71}) has the following form:
\begin{equation}\label{s72}
3\cdot 4=4+8=2\cdot 2+ (0+1+2+(1\cdot1+1\cdot 2+ 1\cdot 2)).
\end{equation}

The Hasse diagram for the theory $T_1\sqcup T_2$ is shown in
Fig.~\ref{fig3}.

Replacing, respectively, $1$ and $2$ limit models of the theories
$T_1$ and $T_2$ by $k>0$ and $m>0$ the equation (\ref{s72}) is
transformed to the following:
$$
(k+2)(m+2)=2\cdot 2+(0+k+m+(k+m+km)).
$$
\end{example}

\begin{example}\label{exx2}\rm
Consider the disjoint union of theory $T_1\sqcup T_2$ in Example
\ref{exx1} and of theory $T_3$ with the first Hasse diagram shown
in Fig.~\ref{fig2}. By Theorem \ref{thchar2} we have the theory
$T_1\sqcup T_2\sqcup T_3$ with $3\cdot 4\cdot 5=60$ countable
models, having the Boolean lattice with $2\cdot 2\cdot 2=8$ prime
models over finite sets and with 52 limit models. The
decomposition formula (\ref{s71}) has the following form:
\begin{equation}\label{s73}
3\cdot 4\cdot 5=8+52=2\cdot 2\cdot 2+ (0+1+2+3+5+7+11+23).
\end{equation}

The Hasse diagram for the theory $T_1\sqcup T_2\sqcup T_3$ is
shown in Fig.~\ref{fig4}.

Replacing, respectively, $1$, $2$, $3$ limit models of the
theories $T_1$, $T_2$, $T_3$ by $k>0$, $m>0$, $n>0$ the equation
(\ref{s73}) is transformed to the following:
$$
(k+2)(m+2)(n+2)=2\cdot 2\cdot
2+(0+k+m+n+(k+m+km)+(k+n+kn)+(m+n+mn)+
$$
$$
+(k+m+n+km+kn+mn+kmn)).
$$
\end{example}

\begin{figure}[t]
\begin{center}
\unitlength 14mm
\begin{picture}(5,1)(-1.14,-0.1)
{\footnotesize \put(0,1){\makebox(0,0)[cc]{$\bullet$}}
\put(0,1){\circle{0.2}} \put(1,0){\circle{0.2}}
\put(1,0){\makebox(0,0)[cc]{$\bullet$}}
\put(2,1){\makebox(0,0)[cc]{$\bullet$}}
\put(2,1){\circle{0.2}}\put(1,2){\circle{0.2}}
\put(1,2){\makebox(0,0)[cc]{$\bullet$}} \put(1,0){\line(1,1){1}}
\put(1,0){\line(-1,1){1}} \put(1,2){\line(-1,-1){1}}
\put(1,2){\line(1,-1){1}} \put(1,-0.2){\makebox(0,0)[tc]{0}}
\put(1,2.2){\makebox(0,0)[bc]{$5$}}
\put(0.06,1.3){\makebox(0,0)[cr]{$1$}}
\put(1.95,1.3){\makebox(0,0)[cl]{$2$}} }
\end{picture}
\hfill \unitlength 14mm
\begin{picture}(5,3.5)(1.35,0.2)
{\footnotesize\put(3,1){\makebox(0,0)[cc]{$\bullet$}}\put(3,1){\circle{0.2}}
\put(3,2){\makebox(0,0)[cc]{$\bullet$}}\put(3,2){\circle{0.2}}
\put(4,0){\makebox(0,0)[cc]{$\bullet$}}\put(4,0){\circle{0.2}}
\put(4,1){\makebox(0,0)[cc]{$\bullet$}}\put(4,1){\circle{0.2}}
\put(4,2){\makebox(0,0)[cc]{$\bullet$}}\put(4,2){\circle{0.2}}
\put(4,3){\makebox(0,0)[cc]{$\bullet$}}\put(4,3){\circle{0.2}}
\put(5,1){\makebox(0,0)[cc]{$\bullet$}}\put(5,1){\circle{0.2}}
\put(5,2){\makebox(0,0)[cc]{$\bullet$}}\put(5,2){\circle{0.2}}
\put(3,1){\line(0,1){1}} \put(5,1){\line(0,1){1}}
\put(4,0){\line(0,1){1}} \put(4,2){\line(0,1){1}}
\put(4,0){\line(1,1){1}} \put(4,1){\line(1,1){1}}
\put(4,0){\line(-1,1){1}} \put(4,1){\line(-1,1){1}}
\put(4,3){\line(-1,-1){1}} \put(4,2){\line(-1,-1){1}}
\put(4,3){\line(1,-1){1}} \put(4,2){\line(1,-1){1}}
\put(4.2,0){\makebox(0,0)[cl]{0}}
\put(3.2,1){\makebox(0,0)[cl]{1}}
\put(4.2,1){\makebox(0,0)[cl]{2}}
\put(5.2,1){\makebox(0,0)[cl]{3}}
\put(3.2,2){\makebox(0,0)[cl]{5}}
\put(4.2,2){\makebox(0,0)[cl]{7}}
\put(5.2,2){\makebox(0,0)[cl]{11}}
\put(4.2,3.1){\makebox(0,0)[cl]{23}}}
\end{picture}
\end{center}
\parbox[t]{0.4\textwidth}{\caption{}\label{fig3}}
\hfill
\parbox[t]{0.4\textwidth}{\caption{}\label{fig4}}

\end{figure}

\begin{example}\label{exx3}\rm
Consider the disjoint union of theory $T_1$ with the third Hasse
diagram shown in Fig.~\ref{fig1}, {\it b} and the theory $T_2$
with the second Hasse diagram shown in Fig.~\ref{fig2}. By Theorem
\ref{thchar2} we have the theory $T_1\sqcup T_2$ with $4\cdot
5=20$ countable models, having the lattice with $3\cdot 3=9$ prime
models over finite sets and with $11$ limit models. The
decomposition formula (\ref{s71}) has the following form:
$$
4\cdot 5=9+11=3\cdot 3+ (0+1+1+1+1+1+3+3).
$$

The Hasse diagram for the theory $T_1\sqcup T_2$ is shown in
Fig.~\ref{fig5}.
\end{example}

\begin{example}\label{exx4}\rm
Consider the disjoint union of theory $T_1\sqcup T_2$ in
Example~\ref{exx3} and the theory $T_3$ with the third Hasse
diagram shown in Fig.~\ref{fig2}. By Theorem \ref{thchar2} we have
the theory $T_1\sqcup T_2\sqcup T_3$ with $4\cdot 5\cdot 5=100$
countable models, having the lattice with $3\cdot 3\cdot 3=27$
prime models over finite sets and with 73 limit models. The
decomposition formula (\ref{s71}) has the following form:
$$
4\cdot 5\cdot 5=27+73=3\cdot 3\cdot 3+(1\cdot 10+2\cdot 2 +3\cdot
4+5\cdot 5+11\cdot 2).
$$

The Hasse diagram for the theory $T_1\sqcup T_2\sqcup T_3$ is
shown in Fig.~\ref{fig6}.
\end{example}

\begin{figure}
\begin{center}
\unitlength 18mm
\begin{picture}(3,1)(-0.69,-0.1)
{\footnotesize \put(0,1){\makebox(0,0)[cc]{$\bullet$}}
\put(0,1){\circle{0.16}} \put(1,0){\circle{0.16}}
\put(1,0){\makebox(0,0)[cc]{$\bullet$}}
\put(2,1){\makebox(0,0)[cc]{$\bullet$}}
\put(2,1){\circle{0.16}}\put(1,2){\circle{0.16}}
\put(1,2){\makebox(0,0)[cc]{$\bullet$}}
\put(1,0){\line(1,1){1}}\put(1.5,0.5){\line(-1,1){1}}
\put(0.5,0.5){\line(1,1){1}} \put(1,0){\line(-1,1){1}}
\put(1,2){\line(-1,-1){1}} \put(1,2){\line(1,-1){1}}
\put(1.5,0.5){\makebox(0,0)[cc]{$\bullet$}}
\put(0.5,0.5){\makebox(0,0)[cc]{$\bullet$}}
\put(1,1){\makebox(0,0)[cc]{$\bullet$}}
\put(0.5,1.5){\makebox(0,0)[cc]{$\bullet$}}
\put(1.5,1.5){\makebox(0,0)[cc]{$\bullet$}}
\put(1.5,0.5){\circle{0.16}} \put(0.5,0.5){\circle{0.16}}
\put(1,1){\circle{0.16}} \put(0.5,1.5){\circle{0.16}}
\put(1.5,1.5){\circle{0.16}} \put(1,-0.17){\makebox(0,0)[tc]{0}}
\put(1,2.15){\makebox(0,0)[bc]{$3$}}
\put(0.06,1.22){\makebox(0,0)[cr]{$1$}}
\put(0.56,0.7){\makebox(0,0)[cr]{$0$}}
\put(1.56,0.7){\makebox(0,0)[cr]{$1$}}
\put(1.56,1.7){\makebox(0,0)[cr]{$1$}}
\put(0.56,1.7){\makebox(0,0)[cr]{$3$}}
\put(1.06,1.22){\makebox(0,0)[cr]{$1$}}
\put(1.95,1.22){\makebox(0,0)[cl]{$1$}} }
\end{picture}
\hfill \unitlength 20mm
\begin{picture}(4,3.5)(1.6,0.2)
{\footnotesize\put(3,1){\makebox(0,0)[cc]{$\bullet$}}\put(3,1){\circle{0.14}}
\put(3,2){\makebox(0,0)[cc]{$\bullet$}}\put(3,2){\circle{0.14}}
\put(4,0){\makebox(0,0)[cc]{$\bullet$}}\put(4,0){\circle{0.14}}
\put(4,1){\makebox(0,0)[cc]{$\bullet$}}\put(4,1){\circle{0.14}}
\put(4.006,2){\makebox(0,0)[cc]{$\bullet$}}\put(4.006,2){\circle{0.08}}
\put(4,3){\makebox(0,0)[cc]{$\bullet$}}\put(4,3){\circle{0.14}}
\put(5,1){\makebox(0,0)[cc]{$\bullet$}}\put(5,1){\circle{0.14}}
\put(5,2){\makebox(0,0)[cc]{$\bullet$}}\put(5,2){\circle{0.14}}
\put(3,1){\line(0,1){1}} \put(5,1){\line(0,1){1}}
\put(4,0){\line(0,1){1}} \put(4,2){\line(0,1){1}}
\put(4,0){\line(1,1){1}} \put(4,1){\line(1,1){1}}
\put(4,0){\line(-1,1){1}} \put(4,1){\line(-1,1){1}}
\put(4,3){\line(-1,-1){1}} \put(4,2){\line(-1,-1){1}}
\put(4,3){\line(1,-1){1}} \put(4,2){\line(1,-1){1}}
\put(4,2.5){\makebox(0,0)[cc]{$\bullet$}}\put(4,2.5){\circle{0.14}}
\put(3.35,0.65){\makebox(0,0)[cc]{$\bullet$}}\put(3.35,0.65){\circle{0.14}}
\put(3.35,1.15){\makebox(0,0)[cc]{$\bullet$}}\put(3.35,1.15){\circle{0.14}}
\put(3.35,1.65){\makebox(0,0)[cc]{$\bullet$}}\put(3.35,1.65){\circle{0.14}}
\put(3.95,1.75){\makebox(0,0)[cc]{$\bullet$}}\put(3.95,1.75){\circle{0.14}}
\put(3.95,1.25){\makebox(0,0)[cc]{$\bullet$}}\put(3.95,1.25){\circle{0.14}}
\put(3.945,2.25){\makebox(0,0)[cc]{$\bullet$}}\put(3.945,2.25){\circle{0.08}}
\put(4.35,1.65){\makebox(0,0)[cc]{$\bullet$}}\put(4.35,1.65){\circle{0.14}}
\put(4.35,2.15){\makebox(0,0)[cc]{$\bullet$}}\put(4.35,2.15){\circle{0.14}}
\put(4.35,2.65){\makebox(0,0)[cc]{$\bullet$}}\put(4.35,2.65){\circle{0.14}}
\put(3.6,1.6){\makebox(0,0)[cc]{$\bullet$}}\put(3.6,1.6){\circle{0.14}}
\put(3.6,2.1){\makebox(0,0)[cc]{$\bullet$}}\put(3.6,2.1){\circle{0.14}}
\put(3.6,2.6){\makebox(0,0)[cc]{$\bullet$}}\put(3.6,2.6){\circle{0.14}}
\put(4.6,0.6){\makebox(0,0)[cc]{$\bullet$}}\put(4.6,0.6){\circle{0.14}}
\put(4.6,1.1){\makebox(0,0)[cc]{$\bullet$}}\put(4.6,1.1){\circle{0.14}}
\put(4.6,1.6){\makebox(0,0)[cc]{$\bullet$}}\put(4.6,1.6){\circle{0.14}}
\put(4,0.5){\makebox(0,0)[cc]{$\bullet$}}\put(4,0.5){\circle{0.14}}
\put(5,1.5){\makebox(0,0)[cc]{$\bullet$}}\put(5,1.5){\circle{0.14}}
\put(3,1.5){\makebox(0,0)[cc]{$\bullet$}}\put(3,1.5){\circle{0.14}}

\put(3.35,0.65){\line(0,1){1}} \put(4.6,0.6){\line(0,1){1}}
\put(3.35,0.65){\line(1,1){1}} \put(4.6,0.6){\line(-1,1){1}}
\put(4,0.5){\line(1,1){1}} \put(4,0.5){\line(-1,1){1}}
\put(3.35,1.65){\line(1,1){1}} \put(4.6,1.6){\line(-1,1){1}}
\put(3,1.5){\line(1,1){1}} \put(5,1.5){\line(-1,1){1}}
\put(3.6,1.6){\line(0,1){1}} \put(4.35,1.65){\line(0,1){1}}
\put(3.35,1.15){\line(1,1){1}} \put(4.6,1.1){\line(-1,1){1}}
\put(3.95,1.25){\line(0,1){1}}

\put(4.1,0){\makebox(0,0)[cl]{0}}
\put(4.1,0.5){\makebox(0,0)[cl]{0}}
\put(4.69,0.6){\makebox(0,0)[cl]{1}}
\put(4.69,1.1){\makebox(0,0)[cl]{1}}
\put(4.69,1.6){\makebox(0,0)[cl]{5}}
\put(3.17,0.65){\makebox(0,0)[cl]{0}}
\put(3.193,1.14){\makebox(0,0)[cl]{0}}
\put(3.193,1.64){\makebox(0,0)[cl]{2}}
\put(3.445,1.63){\makebox(0,0)[cl]{3}}
\put(4.185,1.66){\makebox(0,0)[cl]{1}}

\put(2.83,1){\makebox(0,0)[cl]{1}}
\put(2.83,1.5){\makebox(0,0)[cl]{1}}
\put(2.83,2){\makebox(0,0)[cl]{5}}
\put(4.076,1){\makebox(0,0)[cl]{2}}
\put(5.1,1){\makebox(0,0)[cl]{1}}
\put(5.1,1.5){\makebox(0,0)[cl]{1}}
\put(5.1,2){\makebox(0,0)[cl]{5}}
\put(4.04,1.25){\makebox(0,0)[cl]{1}}
\put(4.04,1.75){\makebox(0,0)[cl]{1}}
\put(3.87,2.0){\makebox(0,0)[cl]{{\tiny 3}}}
\put(3.44,2.1){\makebox(0,0)[cl]{3}}
\put(3.36,2.6){\makebox(0,0)[cl]{11}}
\put(3.825,2.52){\makebox(0,0)[cl]{3}}
\put(4.01,2.25){\makebox(0,0)[cl]{{\tiny 5}}}
\put(4.43,2.15){\makebox(0,0)[cl]{1}}
\put(4.44,2.65){\makebox(0,0)[cl]{5}}
\put(3.929,3.16){\makebox(0,0)[cl]{11}}

}
\end{picture}
\end{center}
\parbox[t]{0.4\textwidth}{\caption{}\label{fig5}}
\hfill
\parbox[t]{0.4\textwidth}{\caption{}\label{fig6}}

\end{figure}

\begin{example}\label{exx5}\rm
Consider the disjoint union of theory $T_1$ with the first Hasse
diagram shown in Fig.~\ref{fig1}, {\it b} and the theory $T_2$
with the second Hasse diagram shown in Fig.~\ref{fig2}. By Theorem
\ref{thchar2} we have the theory $T_1\sqcup T_2$ with $4\cdot
5=20$ countable models, having the lattice with $2\cdot 3=6$ prime
models over finite sets and with $14$ limit models. The
decomposition formula (\ref{s71}) has the following form:
$$
4\cdot 5=6+14=2\cdot 3+ (0+1+1+2+5+5).
$$

The Hasse diagram for the theory $T_1\sqcup T_2$ is shown in
Fig.~\ref{fig7}.
\end{example}

\begin{example}\label{exx6}\rm
Consider the disjoint union of theory $T_1\sqcup T_2$ in
Example~\ref{exx5} and of the theory $T_3$ with the first Hasse
diagram shown in Fig.~\ref{fig2}. By Theorem \ref{thchar2} we have
the theory $T_1\sqcup T_2\sqcup T_3$ with $4\cdot 5\cdot 5=100$
countable models, having the lattice with $2\cdot 2\cdot 3=12$
prime models over finite sets and with 88 limit models. The
decomposition formula (\ref{s71}) has the following form:
$$
4\cdot 5\cdot 5=12+88=2\cdot 2\cdot 3+(1\cdot 2+2\cdot 1 +3\cdot
1+5\cdot 2+7\cdot 2 +11\cdot 1+23\cdot 2).
$$

The Hasse diagram for the theory $T_1\sqcup T_2\sqcup T_3$ is
shown in Fig.~\ref{fig8}.
\end{example}

\begin{figure}[t]
\begin{center}
\unitlength 24mm
\begin{picture}(3,1)(-0.05,-0.1)
{\footnotesize  \put(1,0){\circle{0.16}}
\put(1,0){\makebox(0,0)[cc]{$\bullet$}}
\put(2,1){\makebox(0,0)[cc]{$\bullet$}} \put(2,1){\circle{0.16}}
\put(1,0){\line(1,1){1}}\put(1.5,0.5){\line(-1,1){0.5}}
\put(0.5,0.5){\line(1,1){1}} \put(1,0){\line(-1,1){0.5}}
\put(1.5,1.5){\line(1,-1){0.5}}
\put(1.5,0.5){\makebox(0,0)[cc]{$\bullet$}}
\put(0.5,0.5){\makebox(0,0)[cc]{$\bullet$}}
\put(1,1){\makebox(0,0)[cc]{$\bullet$}}

\put(1.5,1.5){\makebox(0,0)[cc]{$\bullet$}}
\put(1.5,0.5){\circle{0.16}} \put(0.5,0.5){\circle{0.16}}
\put(1,1){\circle{0.16}} \put(1.5,1.5){\circle{0.16}}
\put(1,-0.17){\makebox(0,0)[tc]{0}}

\put(0.56,0.7){\makebox(0,0)[cr]{$2$}}
\put(1.56,0.7){\makebox(0,0)[cr]{$1$}}
\put(1.56,1.7){\makebox(0,0)[cr]{$5$}}

\put(1.06,1.22){\makebox(0,0)[cr]{$5$}}
\put(1.95,1.22){\makebox(0,0)[cl]{$1$}} }
\end{picture}
\hfill \unitlength 24mm
\begin{picture}(3,2.5)(2.65,0.1)
{\footnotesize

\put(4,0){\makebox(0,0)[cc]{$\bullet$}}\put(4,0){\circle{0.14}}
\put(4,1){\makebox(0,0)[cc]{$\bullet$}}\put(4,1){\circle{0.14}}

\put(5,1){\makebox(0,0)[cc]{$\bullet$}}\put(5,1){\circle{0.14}}

\put(5,1){\line(0,1){0.5}} \put(4,0){\line(0,1){0.5}}
 \put(4,0){\line(1,1){1}}
\put(4,0){\line(-1,1){0.5}}

 \put(4.5,1.5){\line(1,-1){0.5}}

\put(3.5,0.5){\makebox(0,0)[cc]{$\bullet$}}\put(3.5,0.5){\circle{0.14}}
\put(3.5,1.0){\makebox(0,0)[cc]{$\bullet$}}\put(3.5,1.0){\circle{0.14}}

\put(4.0,1.5){\makebox(0,0)[cc]{$\bullet$}}\put(4.0,1.5){\circle{0.14}}
\put(4.5,1.5){\makebox(0,0)[cc]{$\bullet$}}\put(4.5,1.5){\circle{0.14}}
\put(4.5,2.0){\makebox(0,0)[cc]{$\bullet$}}\put(4.5,2.0){\circle{0.14}}

\put(4.5,0.5){\makebox(0,0)[cc]{$\bullet$}}\put(4.5,0.5){\circle{0.14}}
\put(4.5,1){\makebox(0,0)[cc]{$\bullet$}}\put(4.5,1){\circle{0.14}}

\put(4,0.5){\makebox(0,0)[cc]{$\bullet$}}\put(4,0.5){\circle{0.14}}
\put(5,1.5){\makebox(0,0)[cc]{$\bullet$}}\put(5,1.5){\circle{0.14}}

\put(3.5,0.5){\line(0,1){0.5}} \put(4.5,0.5){\line(0,1){0.5}}
\put(3.5,0.5){\line(1,1){1}}
\put(4.5,0.5){\line(-1,1){0.5}}\put(4.5,1){\line(-1,1){0.5}}
\put(4,0.5){\line(1,1){1}} \put(4,0.5){\line(-1,1){0.5}}

\put(5,1.5){\line(-1,1){0.5}} \put(4.5,1.5){\line(0,1){0.5}}
\put(3.5,1.0){\line(1,1){1}} \put(4,1){\line(0,1){0.5}}

\put(3.83,0){\makebox(0,0)[cl]{0}}
\put(3.83,0.5){\makebox(0,0)[cl]{3}}
\put(4.61,0.5){\makebox(0,0)[cl]{1}}
\put(4.61,1.0){\makebox(0,0)[cl]{7}}
\put(4.61,1.5){\makebox(0,0)[cl]{5}}
\put(4.45,2.17){\makebox(0,0)[cl]{23}}
\put(3.33,0.5){\makebox(0,0)[cl]{2}}
\put(3.29,1.0){\makebox(0,0)[cl]{11}}

\put(3.83,1){\makebox(0,0)[cl]{5}}
\put(3.78,1.5){\makebox(0,0)[cl]{23}}
\put(5.1,1){\makebox(0,0)[cl]{1}}
\put(5.1,1.5){\makebox(0,0)[cl]{7}}

}
\end{picture}
\end{center}
\parbox[t]{0.4\textwidth}{\caption{}\label{fig7}}
\hfill
\parbox[t]{0.4\textwidth}{\caption{}\label{fig8}}

\end{figure}

\begin{example}\label{exx7}\rm
Consider the disjoint union of theory $T_1$ with the first Hasse
diagram shown in Fig.~\ref{fig1}, {\it b} and of the theory $T_2$
with the third Hasse diagram shown in Fig.~\ref{fig2}. By Theorem
\ref{thchar2} we have the theory $T_1\sqcup T_2$ with $4\cdot
5=20$ countable models, having the lattice with $2\cdot 3=6$ prime
models over finite sets and with $14$ limit models. The
decomposition formula (\ref{s71}) has the following form:
$$
4\cdot 5=6+14=2\cdot 3+ (2\cdot 3+8).
$$

The Hasse diagram for the theory $T_1\sqcup T_2$ is shown in
Fig.~\ref{fig9}.
\end{example}

\begin{example}\label{exx8}\rm
Consider the disjoint union of theory $T_1\sqcup T_2$ in
Example~\ref{exx7} and of the theory $T_3$ with the first Hasse
diagram shown in Fig.~\ref{fig2}. By Theorem \ref{thchar2} we have
the theory $T_1\sqcup T_2\sqcup T_3$ with $4\cdot 5\cdot 5=100$
countable models, having the lattice with $2\cdot 2\cdot 3=12$
prime models over finite sets and with 88 limit models. The
decomposition formula (\ref{s71}) has the following form:
$$
4\cdot 5\cdot 5=12+88=2\cdot 2\cdot 3+(2\cdot 3 +3\cdot
2+8+11\cdot 3+35).
$$

The Hasse diagram for the theory $T_1\sqcup T_2\sqcup T_3$ is
shown in Fig.~\ref{fig10}.
\end{example}

\begin{figure}[t]
\begin{center}
\unitlength 24mm
\begin{picture}(3,1)(-0.05,-0.1)
{\footnotesize  \put(1,0){\circle{0.16}}
\put(1,0){\makebox(0,0)[cc]{$\bullet$}}
\put(2,1){\makebox(0,0)[cc]{$\bullet$}} \put(2,1){\circle{0.16}}
\put(1,0){\line(1,1){1}}\put(1.5,0.5){\line(-1,1){0.5}}
\put(0.5,0.5){\line(1,1){1}} \put(1,0){\line(-1,1){0.5}}
\put(1.5,1.5){\line(1,-1){0.5}}
\put(1.5,0.5){\makebox(0,0)[cc]{$\bullet$}}
\put(0.5,0.5){\makebox(0,0)[cc]{$\bullet$}}
\put(1,1){\makebox(0,0)[cc]{$\bullet$}}

\put(1.5,1.5){\makebox(0,0)[cc]{$\bullet$}}
\put(1.5,0.5){\circle{0.16}} \put(0.5,0.5){\circle{0.16}}
\put(1,1){\circle{0.16}} \put(1.5,1.5){\circle{0.16}}
\put(1,-0.17){\makebox(0,0)[tc]{0}}

\put(0.56,0.7){\makebox(0,0)[cr]{$2$}}
\put(1.56,0.7){\makebox(0,0)[cr]{$0$}}
\put(1.56,1.7){\makebox(0,0)[cr]{$8$}}

\put(1.06,1.22){\makebox(0,0)[cr]{$2$}}
\put(1.95,1.22){\makebox(0,0)[cl]{$2$}} }
\end{picture}
\hfill \unitlength 24mm
\begin{picture}(3,2.5)(2.65,0.1)
{\footnotesize

\put(4,0){\makebox(0,0)[cc]{$\bullet$}}\put(4,0){\circle{0.14}}
\put(4,1){\makebox(0,0)[cc]{$\bullet$}}\put(4,1){\circle{0.14}}

\put(5,1){\makebox(0,0)[cc]{$\bullet$}}\put(5,1){\circle{0.14}}

\put(5,1){\line(0,1){0.5}} \put(4,0){\line(0,1){0.5}}
 \put(4,0){\line(1,1){1}}
\put(4,0){\line(-1,1){0.5}}

 \put(4.5,1.5){\line(1,-1){0.5}}

\put(3.5,0.5){\makebox(0,0)[cc]{$\bullet$}}\put(3.5,0.5){\circle{0.14}}
\put(3.5,1.0){\makebox(0,0)[cc]{$\bullet$}}\put(3.5,1.0){\circle{0.14}}

\put(4.0,1.5){\makebox(0,0)[cc]{$\bullet$}}\put(4.0,1.5){\circle{0.14}}
\put(4.5,1.5){\makebox(0,0)[cc]{$\bullet$}}\put(4.5,1.5){\circle{0.14}}
\put(4.5,2.0){\makebox(0,0)[cc]{$\bullet$}}\put(4.5,2.0){\circle{0.14}}

\put(4.5,0.5){\makebox(0,0)[cc]{$\bullet$}}\put(4.5,0.5){\circle{0.14}}
\put(4.5,1){\makebox(0,0)[cc]{$\bullet$}}\put(4.5,1){\circle{0.14}}

\put(4,0.5){\makebox(0,0)[cc]{$\bullet$}}\put(4,0.5){\circle{0.14}}
\put(5,1.5){\makebox(0,0)[cc]{$\bullet$}}\put(5,1.5){\circle{0.14}}

\put(3.5,0.5){\line(0,1){0.5}} \put(4.5,0.5){\line(0,1){0.5}}
\put(3.5,0.5){\line(1,1){1}}
\put(4.5,0.5){\line(-1,1){0.5}}\put(4.5,1){\line(-1,1){0.5}}
\put(4,0.5){\line(1,1){1}} \put(4,0.5){\line(-1,1){0.5}}

\put(5,1.5){\line(-1,1){0.5}} \put(4.5,1.5){\line(0,1){0.5}}
\put(3.5,1.0){\line(1,1){1}} \put(4,1){\line(0,1){0.5}}

\put(3.83,0){\makebox(0,0)[cl]{0}}
\put(3.83,0.5){\makebox(0,0)[cl]{3}}
\put(4.61,0.5){\makebox(0,0)[cl]{0}}
\put(4.61,1.0){\makebox(0,0)[cl]{3}}
\put(4.61,1.5){\makebox(0,0)[cl]{8}}
\put(4.45,2.17){\makebox(0,0)[cl]{35}}
\put(3.33,0.5){\makebox(0,0)[cl]{2}}
\put(3.29,1.0){\makebox(0,0)[cl]{11}}

\put(3.83,1){\makebox(0,0)[cl]{2}}
\put(3.78,1.5){\makebox(0,0)[cl]{11}}
\put(5.1,1){\makebox(0,0)[cl]{2}}
\put(5.1,1.5){\makebox(0,0)[cl]{11}}

}
\end{picture}
\end{center}
\parbox[t]{0.4\textwidth}{\caption{}\label{fig9}}
\hfill
\parbox[t]{0.4\textwidth}{\caption{}\label{fig10}}

\end{figure}

\begin{example}\label{exx9}\rm
Consider the disjoint union of theory $T_1\sqcup T_2$ in
Example~\ref{exx3} and of the theory $T_3$ with the first Hasse
diagram shown in Fig.~\ref{fig1}, {\it b}. By Theorem
\ref{thchar2} we have the theory $T_1\sqcup T_2\sqcup T_3$ with
$4\cdot 4\cdot 5=80$ countable models, having the lattice with
$3\cdot 3\cdot 2=18$ prime models over finite sets and with $62$
limit models. The decomposition formula (\ref{s71}) has the
following form:
$$
4\cdot 4\cdot 5=18+62=3\cdot 3\cdot 2+(1\cdot 5+2\cdot 2 +3\cdot
2+5\cdot 5+11\cdot 2).
$$

The Hasse diagram for the theory $T_1\sqcup T_2\sqcup T_3$ is
shown in Fig.~\ref{fig11}.
\end{example}

\begin{example}\label{exx10}\rm
Consider the disjoint union $T_1\sqcup T_2\sqcup T_3$ of theory
$T_1$ with the third Hasse diagram shown in Fig.~\ref{fig1}, {\it
b}, of theory $T_2$ with the first Hasse diagram shown in
Fig.~\ref{fig2}, and of theory $T_3$ with the third Hasse diagram
shown in Fig.~\ref{fig2}. By Theorem \ref{thchar2} we have the
theory $T_1\sqcup T_2\sqcup T_3$ with $4\cdot 5\cdot 5=100$
countable models, having the lattice with $3\cdot 3\cdot 2=18$
prime models over finite sets and with $82$ limit models. The
decomposition formula (\ref{s71}) has the following form:
$$
4\cdot 4\cdot 5=18+62=3\cdot 3\cdot 2+(1\cdot 2+2\cdot 2 +3\cdot
4+5\cdot 1+7\cdot 2+11\cdot 2+23\cdot 1).
$$

The Hasse diagram for the theory $T_1\sqcup T_2\sqcup T_3$ is
shown in Fig.~\ref{fig12}.
\end{example}

\begin{figure}[t]
\begin{center}
\unitlength 20mm
\begin{picture}(3.5,2.95)(2.52,0.0)
{\footnotesize\put(3,1){\makebox(0,0)[cc]{$\bullet$}}\put(3,1){\circle{0.14}}

\put(4,0){\makebox(0,0)[cc]{$\bullet$}}\put(4,0){\circle{0.14}}
\put(4,1){\makebox(0,0)[cc]{$\bullet$}}\put(4,1){\circle{0.14}}
\put(4.0,2){\makebox(0,0)[cc]{$\bullet$}}\put(4.0,2){\circle{0.14}}

\put(5,1){\makebox(0,0)[cc]{$\bullet$}}\put(5,1){\circle{0.14}}

\put(3,1){\line(0,1){0.5}} \put(5,1){\line(0,1){0.5}}
\put(4,0){\line(0,1){0.5}} \put(4,2){\line(0,1){0.5}}
\put(4,0){\line(1,1){1}} \put(4,0){\line(-1,1){1}}
 \put(4,2){\line(-1,-1){1}}
 \put(4,2){\line(1,-1){1}}
\put(4,2.5){\makebox(0,0)[cc]{$\bullet$}}\put(4,2.5){\circle{0.14}}
\put(3.5,0.5){\makebox(0,0)[cc]{$\bullet$}}\put(3.5,0.5){\circle{0.14}}
\put(3.5,1.0){\makebox(0,0)[cc]{$\bullet$}}\put(3.5,1.0){\circle{0.14}}
\put(3.5,1.5){\makebox(0,0)[cc]{$\bullet$}}\put(3.5,1.5){\circle{0.14}}
\put(3.5,2.0){\makebox(0,0)[cc]{$\bullet$}}\put(3.5,2.0){\circle{0.14}}
\put(4.0,1.5){\makebox(0,0)[cc]{$\bullet$}}\put(4.0,1.5){\circle{0.14}}
\put(4.5,1.5){\makebox(0,0)[cc]{$\bullet$}}\put(4.5,1.5){\circle{0.14}}
\put(4.5,2.0){\makebox(0,0)[cc]{$\bullet$}}\put(4.5,2.0){\circle{0.14}}

\put(4.5,0.5){\makebox(0,0)[cc]{$\bullet$}}\put(4.5,0.5){\circle{0.14}}
\put(4.5,1){\makebox(0,0)[cc]{$\bullet$}}\put(4.5,1){\circle{0.14}}

\put(4,0.5){\makebox(0,0)[cc]{$\bullet$}}\put(4,0.5){\circle{0.14}}
\put(5,1.5){\makebox(0,0)[cc]{$\bullet$}}\put(5,1.5){\circle{0.14}}
\put(3,1.5){\makebox(0,0)[cc]{$\bullet$}}\put(3,1.5){\circle{0.14}}

\put(3.5,0.5){\line(0,1){0.5}} \put(4.5,0.5){\line(0,1){0.5}}
\put(3.5,0.5){\line(1,1){1}}
\put(4.5,0.5){\line(-1,1){1}}\put(4.5,1){\line(-1,1){1}}
\put(4,0.5){\line(1,1){1}} \put(4,0.5){\line(-1,1){1}}
 \put(3,1.5){\line(1,1){1}}
\put(5,1.5){\line(-1,1){1}} \put(3.5,1.5){\line(0,1){0.5}}
\put(4.5,1.5){\line(0,1){0.5}} \put(3.5,1.0){\line(1,1){1}}
\put(4,1){\line(0,1){0.5}}

\put(4.1,0){\makebox(0,0)[cl]{0}}
\put(4.1,0.5){\makebox(0,0)[cl]{2}}
\put(4.61,0.5){\makebox(0,0)[cl]{0}}
\put(4.61,1.0){\makebox(0,0)[cl]{1}}
\put(4.61,1.5){\makebox(0,0)[cl]{1}}
\put(4.61,2.0){\makebox(0,0)[cl]{5}}
\put(3.33,0.5){\makebox(0,0)[cl]{0}}
\put(3.33,1.0){\makebox(0,0)[cl]{2}}
\put(3.33,1.5){\makebox(0,0)[cl]{3}}
\put(3.33,2.0){\makebox(0,0)[cl]{11}}

\put(2.83,1){\makebox(0,0)[cl]{1}}
\put(2.83,1.5){\makebox(0,0)[cl]{5}}

\put(4.1,1){\makebox(0,0)[cl]{1}}
\put(4.1,1.5){\makebox(0,0)[cl]{5}}
\put(5.1,1){\makebox(0,0)[cl]{1}}
\put(5.1,1.5){\makebox(0,0)[cl]{5}}

\put(4.1,2.0){\makebox(0,0)[cl]{3}}

\put(3.93,2.67){\makebox(0,0)[cl]{11}}

}
\end{picture}
\hfill\unitlength 20mm
\begin{picture}(4,2.95)(1.65,0.0)
{\footnotesize\put(3,1){\makebox(0,0)[cc]{$\bullet$}}\put(3,1){\circle{0.14}}

\put(4,0){\makebox(0,0)[cc]{$\bullet$}}\put(4,0){\circle{0.14}}
\put(4,1){\makebox(0,0)[cc]{$\bullet$}}\put(4,1){\circle{0.14}}
\put(4.0,2){\makebox(0,0)[cc]{$\bullet$}}\put(4.0,2){\circle{0.14}}

\put(5,1){\makebox(0,0)[cc]{$\bullet$}}\put(5,1){\circle{0.14}}

\put(3,1){\line(0,1){0.5}} \put(5,1){\line(0,1){0.5}}
\put(4,0){\line(0,1){0.5}} \put(4,2){\line(0,1){0.5}}
\put(4,0){\line(1,1){1}} \put(4,0){\line(-1,1){1}}
 \put(4,2){\line(-1,-1){1}}
 \put(4,2){\line(1,-1){1}}
\put(4,2.5){\makebox(0,0)[cc]{$\bullet$}}\put(4,2.5){\circle{0.14}}
\put(3.5,0.5){\makebox(0,0)[cc]{$\bullet$}}\put(3.5,0.5){\circle{0.14}}
\put(3.5,1.0){\makebox(0,0)[cc]{$\bullet$}}\put(3.5,1.0){\circle{0.14}}
\put(3.5,1.5){\makebox(0,0)[cc]{$\bullet$}}\put(3.5,1.5){\circle{0.14}}
\put(3.5,2.0){\makebox(0,0)[cc]{$\bullet$}}\put(3.5,2.0){\circle{0.14}}
\put(4.0,1.5){\makebox(0,0)[cc]{$\bullet$}}\put(4.0,1.5){\circle{0.14}}
\put(4.5,1.5){\makebox(0,0)[cc]{$\bullet$}}\put(4.5,1.5){\circle{0.14}}
\put(4.5,2.0){\makebox(0,0)[cc]{$\bullet$}}\put(4.5,2.0){\circle{0.14}}

\put(4.5,0.5){\makebox(0,0)[cc]{$\bullet$}}\put(4.5,0.5){\circle{0.14}}
\put(4.5,1){\makebox(0,0)[cc]{$\bullet$}}\put(4.5,1){\circle{0.14}}

\put(4,0.5){\makebox(0,0)[cc]{$\bullet$}}\put(4,0.5){\circle{0.14}}
\put(5,1.5){\makebox(0,0)[cc]{$\bullet$}}\put(5,1.5){\circle{0.14}}
\put(3,1.5){\makebox(0,0)[cc]{$\bullet$}}\put(3,1.5){\circle{0.14}}

\put(3.5,0.5){\line(0,1){0.5}} \put(4.5,0.5){\line(0,1){0.5}}
\put(3.5,0.5){\line(1,1){1}}
\put(4.5,0.5){\line(-1,1){1}}\put(4.5,1){\line(-1,1){1}}
\put(4,0.5){\line(1,1){1}} \put(4,0.5){\line(-1,1){1}}
 \put(3,1.5){\line(1,1){1}}
\put(5,1.5){\line(-1,1){1}} \put(3.5,1.5){\line(0,1){0.5}}
\put(4.5,1.5){\line(0,1){0.5}} \put(3.5,1.0){\line(1,1){1}}
\put(4,1){\line(0,1){0.5}}

\put(4.1,0){\makebox(0,0)[cl]{0}}
\put(4.1,0.5){\makebox(0,0)[cl]{3}}
\put(4.61,0.5){\makebox(0,0)[cl]{0}}
\put(4.61,1.0){\makebox(0,0)[cl]{3}}
\put(4.61,1.5){\makebox(0,0)[cl]{2}}
\put(4.61,2.0){\makebox(0,0)[cl]{11}}
\put(3.33,0.5){\makebox(0,0)[cl]{0}}
\put(3.33,1.0){\makebox(0,0)[cl]{3}}
\put(3.33,1.5){\makebox(0,0)[cl]{1}}
\put(3.33,2.0){\makebox(0,0)[cl]{7}}

\put(2.83,1){\makebox(0,0)[cl]{1}}
\put(2.83,1.5){\makebox(0,0)[cl]{7}}

\put(4.1,1){\makebox(0,0)[cl]{0}}
\put(4.1,1.5){\makebox(0,0)[cl]{3}}
\put(5.1,1){\makebox(0,0)[cl]{2}}
\put(5.1,1.5){\makebox(0,0)[cl]{11}}

\put(4.1,2.0){\makebox(0,0)[cl]{5}}

\put(3.93,2.67){\makebox(0,0)[cl]{23}}

}
\end{picture}

\end{center}
\parbox[t]{0.4\textwidth}{\caption{}\label{fig11}}
\hfill
\parbox[t]{0.4\textwidth}{\caption{}\label{fig12}}

\end{figure}

\begin{example}\label{exx11}\rm
Consider the disjoint union of theory $T_1$ with the second Hasse
diagram shown in Fig.~\ref{fig1}, {\it b}, where one limit model
is replaced by $k>0$ ones, and of theory $T_2$ with the Hasse
diagram shown in Fig.~\ref{fig1}, {\it a}, where one limit model
is replaced by $m>0$ ones. By Theorem \ref{thchar2} we have the
theory $T_1\sqcup T_2\sqcup T_3$ with $(k+3)(m+2)$ countable
models, having the Hasse diagram with $3\cdot 2=6$ prime models
over finite sets and with $k+m+(k+2m+km)$ limit models. The
decomposition formula (\ref{s71}) has the following form:
$$
(k+3)(m+2)=3\cdot 2+(k+m+(k+2m+km)).
$$

The Hasse diagram for the theory $T_1\sqcup T_2$ is shown in
Fig.~\ref{fig13}.
\end{example}

\begin{example}\label{exx12}\rm
Consider the disjoint union of theory $T_1$ with the second Hasse
diagram shown in Fig.~\ref{fig1}, {\it b}, where one limit model
is replaced by $k>0$ ones, and of theory $T_2$ with similar Hasse
diagram, where one limit model is replaced by $m>0$ ones. By
Theorem \ref{thchar2} we have the theory $T_1\sqcup T_2$ with
$(k+3)(m+3)$ countable models, having the Hasse diagram with
$3\cdot 3=9$ prime models over finite sets and with
$k+m+(2k+2m+km)$ limit models. The decomposition formula
(\ref{s71}) has the following form:
$$
(k+3)(m+3)=3\cdot 3+(k+m+(2k+2m+km)).
$$

The Hasse diagram for the theory $T_1\sqcup T_2$ is shown in
Fig.~\ref{fig14}.
\end{example}

In the latter two examples quotients with respect to $\sim_{\rm
RK}$ produce Boolean lattices with four elements.

\begin{figure}
\vspace{5mm}
\begin{center}
\unitlength 4mm
\begin{picture}(7,5)(25.5,8.5)
{\footnotesize

\put(14.5,2.5){\line(-2,5){2}} \put(14.5,2.5){\line(2,5){2}}
\put(14.5,2.5){\makebox(0,0)[cc]{$\bullet$}}
\put(12.5,7.5){\makebox(0,0)[cc]{$\bullet$}}
\put(16.5,7.5){\makebox(0,0)[cc]{$\bullet$}}
\put(14.5,2.5){\circle{0.6}} \put(14.5,7.5){\oval(5,1)}
\put(15.5,2.5){\makebox(0,0)[cr]{$0$}}
\put(11.7,7.5){\makebox(0,0)[cr]{$k$}}

\put(20.5,6.5){\line(-2,5){2}} \put(20.5,6.5){\line(2,5){2}}
\put(20.5,6.5){\makebox(0,0)[cc]{$\bullet$}}
\put(18.5,11.5){\makebox(0,0)[cc]{$\bullet$}}
\put(22.5,11.5){\makebox(0,0)[cc]{$\bullet$}}
\put(20.5,6.5){\circle{0.6}} \put(20.5,11.5){\oval(5,1)}
\put(21.7,6.5){\makebox(0,0)[cr]{$m$}}
\put(22.9,12.5){\makebox(0,0)[cr]{$k+2m+km$}}

\put(14.5,2.5){\line(3,2){6}} \put(12.5,7.5){\line(3,2){6}}
\put(16.5,7.5){\line(3,2){6}}

\put(12.5,7.5){\line(5,2){10}} \put(16.5,7.5){\line(2,4){2}}

}
\end{picture}

\begin{picture}(7,5)(7.5,3.5)
{\footnotesize

\put(14.5,2.5){\line(-2,5){2}} \put(14.5,2.5){\line(2,5){2}}
\put(14.5,2.5){\makebox(0,0)[cc]{$\bullet$}}
\put(12.5,7.5){\makebox(0,0)[cc]{$\bullet$}}
\put(16.5,7.5){\makebox(0,0)[cc]{$\bullet$}}
\put(14.5,2.5){\circle{0.6}} \put(14.5,7.5){\oval(5,1)}
\put(15.5,2.5){\makebox(0,0)[cr]{$0$}}
\put(11.7,7.5){\makebox(0,0)[cr]{$k$}}

\put(20.5,6.5){\line(-2,5){2}} \put(20.5,6.5){\line(2,5){2}}
\put(20.5,6.5){\makebox(0,0)[cc]{$\bullet$}}
\put(18.5,11.5){\makebox(0,0)[cc]{$\bullet$}}
\put(22.5,11.5){\makebox(0,0)[cc]{$\bullet$}}
%\put(20.5,6.5){\circle{0.6}}
\put(20.5,11.5){\oval(5,1)}
%\put(21.7,6.5){\makebox(0,0)[cr]{$m$}}
%\put(22.9,12.5){\makebox(0,0)[cr]{$k+2m+km$}}

\put(26.5,6.5){\line(-2,5){2}} \put(26.5,6.5){\line(2,5){2}}
\put(26.5,6.5){\makebox(0,0)[cc]{$\bullet$}}
\put(24.5,11.5){\makebox(0,0)[cc]{$\bullet$}}
\put(28.5,11.5){\makebox(0,0)[cc]{$\bullet$}}
%\put(26.5,6.5){\circle{0.6}}
\put(26.5,11.5){\oval(5,1)}
\put(27.85,6.5){\makebox(0,0)[cr]{$m$}}
\put(26.1,12.8){\makebox(0,0)[cr]{$2k+2m+km$}}

\put(26.5,6.5){\line(-2,5){2}} \put(26.5,6.5){\line(2,5){2}}

\put(14.5,2.5){\line(3,2){6}} \put(12.5,7.5){\line(3,2){6}}
\put(16.5,7.5){\line(3,2){6}}

\put(14.5,2.5){\line(6,2){12}}

\put(12.5,7.5){\line(6,2){12}} \put(16.5,7.5){\line(6,2){12}}

\put(12.5,7.5){\line(5,2){10}} \put(16.5,7.5){\line(4,2){8}}

\put(12.5,7.5){\line(4,1){16}}

\put(16.5,7.5){\line(2,4){2}}

\put(23.45,6.5){\oval(7,1)}

\put(23.5,11.5){\oval(12,1.5)}

}
\end{picture}

\end{center}
\parbox[t]{0.4\textwidth}{\caption{}\label{fig13}}
\hfill
\parbox[t]{0.4\textwidth}{\caption{}\label{fig14}}
\end{figure}

\end{document}